\documentclass{article}

\usepackage{color}
\usepackage{graphicx}
\usepackage{sgame}
\usepackage{multirow}

\title{A Short Note on Nonlinear Games on a Grid}
\author{Stewart D. Johnson \thanks{\tt sjohnson@williams.edu}  \\Department of Mathematics and Statistics\\Williams College, Williamstown, MA 01267}

\begin{document}             
\openup 4pt
\maketitle

\begin{abstract}
\noindent Players are arranged on a regular lattice and coded with a specific strategy for a pre-defined game. Each player sums their payoffs from playing the game with each of their neighbors, and then adopts the strategy of the most successful player in the neighborhood. Dynamics are thus determined by the relative ranks of all possible payoff sums.
\medskip

\noindent Linear sums of payoffs, however, generate only a small proportion of all possible rankings. Allowing for any ranking (motivated by Conway's Game of Life) creates a rich array of dynamics.

\bigskip

\end{abstract}

\vfil\break
\renewcommand{\thefootnote}{\fnsymbol{footnote}}

\section{Games on a Grid}

Games on a grid have been extensively studied (see \cite{NowMay}\cite{LinNor}\cite{KillDoe}\cite{NowSig}\cite{Stern}) with interesting dynamics and the surprising result that strictly dominated strategies, like cooperation in Prisoner's Dilemma, can persist in a spatial environment.

The simplest idea is that you start with and underlying $2\times2$ symmetric game
\begin{center}
\begin{game}{2}{2}
& coop & defect\\
coop   & $\begin{array}{c} \qquad a \\ a \qquad  \end{array}$ & $\begin{array}{c} \qquad b  \\ c \qquad  \end{array}$\\
defect & $\begin{array}{c} \qquad c  \\ b \qquad  \end{array}$ & $\begin{array}{c} \qquad d  \\ d \qquad  \end{array}$
\end{game}
\end{center}

and consider an $N\times M$ grid of players. Using wrap-around to make it a torus is convenient to avoid boundary effects.

Dynamics proceed by repeating two steps:

STEP 1) Each player is initially programmed to play either coop or defect, plays the game with each of their 8 neighbors, and sums their own payoffs.

STEP 2) Each player compares their sum to each of their 8 neighbors. If any neighbor has a higher sum using a different strategy, the player will adopt the strategy of this more successful neighbor.

For generic payoffs, there will be no ties between players using different strategies. For non-generic 2x2 games, one can adopt the convention of preferring to keep or to change one's strategy. Introducing a stochastic element to break ties is another solution, but will not be considered here.

Hexagonal grids, up-down-left-right neighborhoods, games on a regular graphs, and games with more than one strategy are easily considered with the same analysis as the following.

\section{Analysis}

Suppose each player has $N$ neighbors, has strategies numbered $0$ or $1$, and the underlying game is:
\begin{center}
\begin{game}{2}{2}
& 0 &  1 \\
0  &  a& b  \\
1 &  c & d
\end{game}
\end{center}

In STEP 1, each player counts the number of neighbors with strategy $1$, obtaining a count $0\le k \le N$. The player's own strategy is either $1$ or $0$, so there are $2(N+1)$ possible net payoffs:

$$\begin{array}{lc}
\hbox{Payoff for player of type 0:\ }& (N-k)a+kb\\
\hbox{Payoff for player of type 1:\ }& (N-k)c+kd
\end{array}
$$

Generically, these will be $2(N+1)$ distinct values (i.e. if $a,b,c,d$ are real numbers selected from any continuous probability distribution, there is a zero probability that two of these expressions will be equal), and hence a player's payoff uniquely determines their strategy.

In STEP 2, a player will compare their net payoff to that of their neighbors. The only thing that matters is the relative ranks of the net payoffs.

For example, consider the standard rectangular grid with eight neighbors and the
prisoner's dilemma:
\begin{center}
\begin{game}{2}{2}
& 0 (coop) & 1 (defect)\\
0 (coop)   &  1.0& 0.1  \\
1 (defect) &  1.9 & 0.3
\end{game}
\end{center}

\vfil\break

The possible payoffs are:

\begin{tabular}{c|c|c|c|c|c|c|c|c|c|}
\multicolumn{1}{c}{} & \multicolumn{9}{c}{\hbox{Number of neighbors of type 1}}\\
\multicolumn{1}{c}{} &\multicolumn{1}{c}{0}&\multicolumn{1}{c}{1}&\multicolumn{1}{c}{2}&\multicolumn{1}{c}{3} &\multicolumn{1}{c}{4}&\multicolumn{1}{c}{5}&\multicolumn{1}{c}{6}&\multicolumn{1}{c}{7}&\multicolumn{1}{c}{8}\\
\cline{2-10}
\hbox{ Player type 0:\ } &8.0&7.1&6.2&5.3&4.4&3.5&2.6&1.7&0.8 \\
\cline{2-10}
\hbox{ Player type 1:\ } &15.2&13.6&12.0&10.4&8.8&7.2&5.6&4.0&2.4 \\
\cline{2-10}
\end{tabular}

The resulting rankings create a rank matrix:

\setlength{\tabcolsep}{8pt}
\begin{tabular}{c|c|c|c|c|c|c|c|c|c|}
\multicolumn{1}{c}{} & \multicolumn{9}{c}{\hbox{Number of neighbors of type 1}}\\
\multicolumn{1}{c}{} &\multicolumn{1}{c}{0}&\multicolumn{1}{c}{1}&\multicolumn{1}{c}{2}&\multicolumn{1}{c}{3} &\multicolumn{1}{c}{4}&\multicolumn{1}{c}{5}&\multicolumn{1}{c}{6}&\multicolumn{1}{c}{7}&\multicolumn{1}{c}{8}\\
\cline{2-10}
\hbox{ Player type 0:\ } &13 &11 &10 &8 &7 &5 &4 &2 &1       \\
\cline{2-10}
\hbox{ Player type 1:\ } &18 &17 &16 &15 &14 &12 &9 &6 &3  \\
\cline{2-10}
\end{tabular}

The ranking determines the game dynamics; once the ranking matrix is determined, the game matrix and net payoffs can be dispensed with.

Note that the net payoff depends linearly on the count of neighbors, and hence both rows of the rank matrix are monotone.

\section{Nonlinearity}

The most well know cellular automata is Conway's Game of Life, whose dynamics also depend on a count of neighbors. However the dependence is not monotone: too many or too few occupied neighbors leads to the demise of an occupied cell.

The core idea of nonlinear spatial games is to consider all possible rank matrices, regardless of monotonicity or even any underlying game.  That is, allow the $2\times (N+1)$ rank matrix to contain any arrangement of the ranks $1,2,\ldots, 2(N+1)$.

\vfil\break

For example, the rank mtx:

\begin{tabular}{c|c|c|c|c|c|c|c|c|c|}
\multicolumn{1}{c}{} & \multicolumn{9}{c}{\hbox{Number of neighbors of type 1}}\\
\multicolumn{1}{c}{} &\multicolumn{1}{c}{0}&\multicolumn{1}{c}{1}&\multicolumn{1}{c}{2}&\multicolumn{1}{c}{3} &\multicolumn{1}{c}{4}&\multicolumn{1}{c}{5}&\multicolumn{1}{c}{6}&\multicolumn{1}{c}{7}&\multicolumn{1}{c}{8}\\
\cline{2-10}
\hbox{ Player type 0:\ } &12&8&16&14&9&3&6&1&10 \\
\cline{2-10}
\hbox{ Player type 1:\ } &2&11&18&17&5&13&4&15&7 \\
\cline{2-10}
\end{tabular}

leads to interesting dynamics:

\begin{center}
\includegraphics[width=1in]{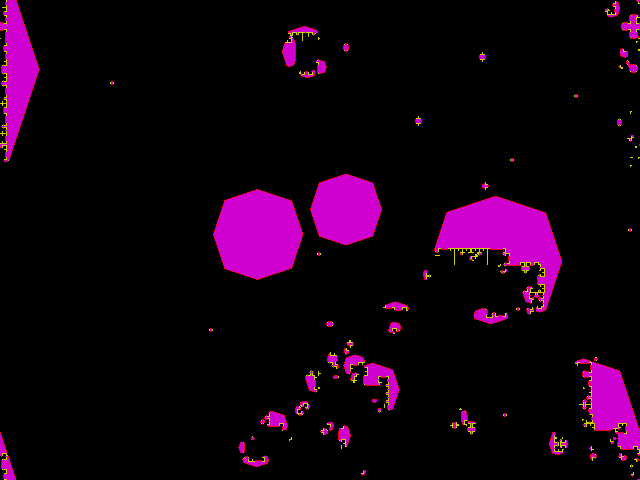}
\includegraphics[width=1in]{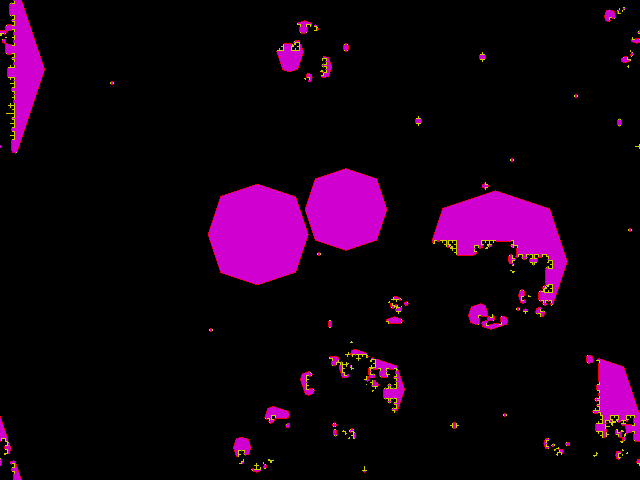}
\includegraphics[width=1in]{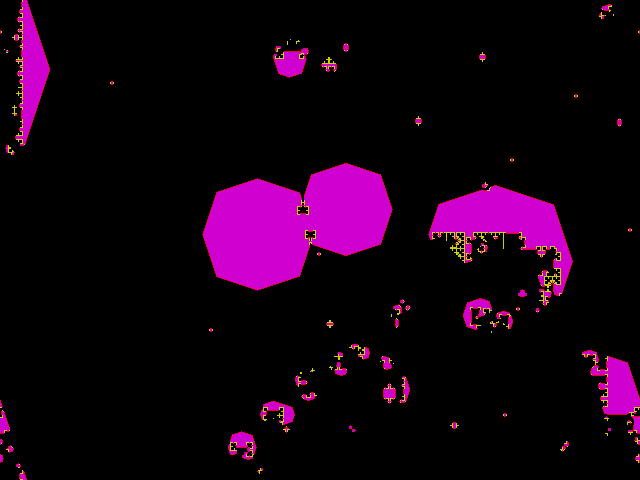}
\includegraphics[width=1in]{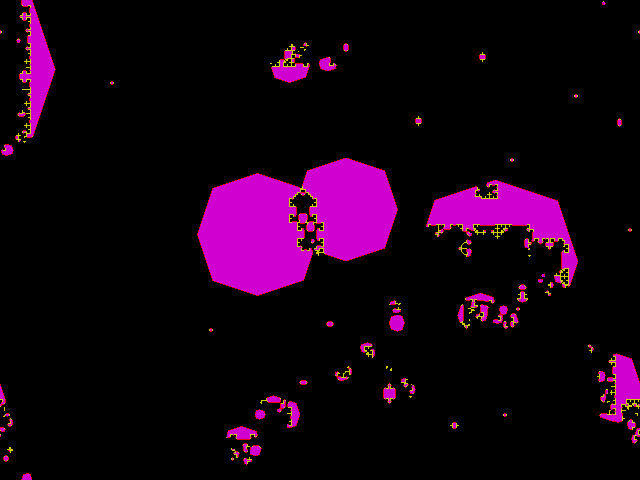}
\end{center}

\begin{center}
\includegraphics[width=1in]{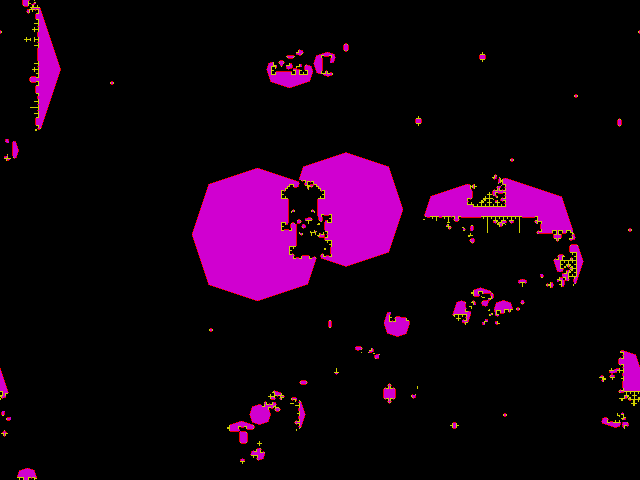}
\includegraphics[width=1in]{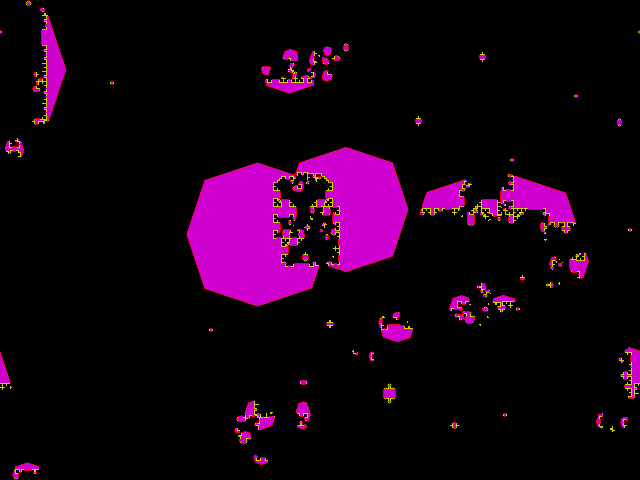}
\includegraphics[width=1in]{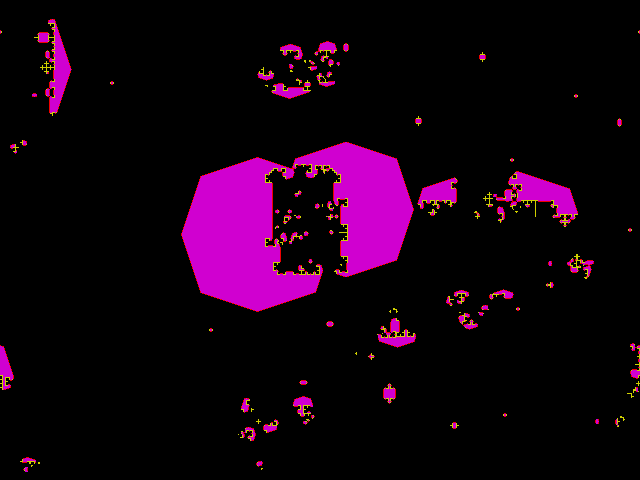}
\includegraphics[width=1in]{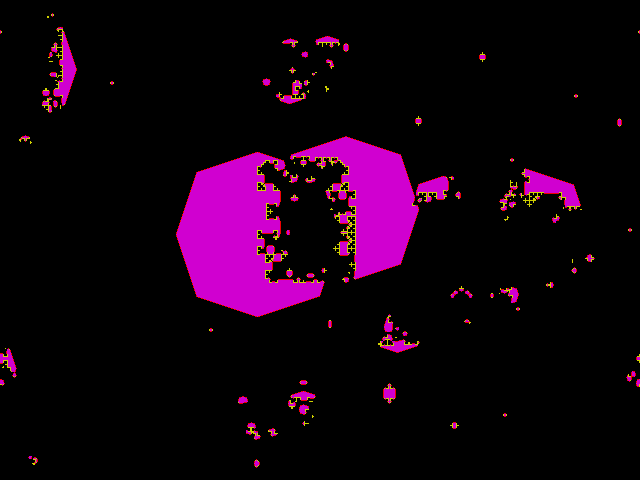}
\end{center}

\begin{center}
\includegraphics[width=1in]{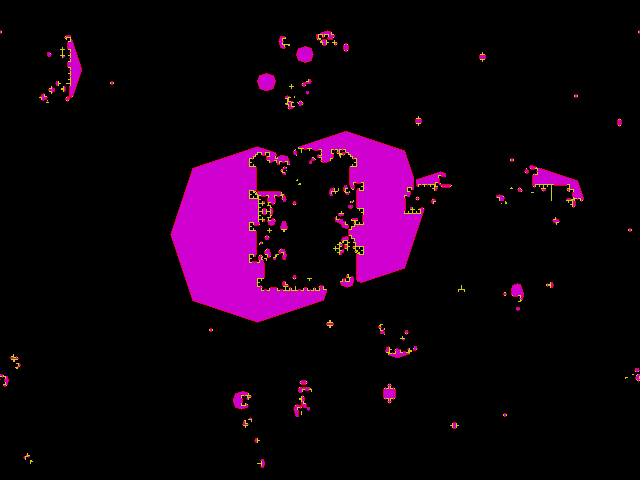}
\includegraphics[width=1in]{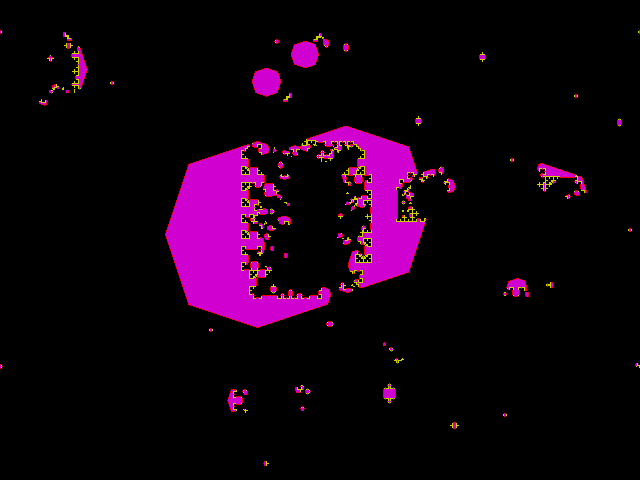}
\includegraphics[width=1in]{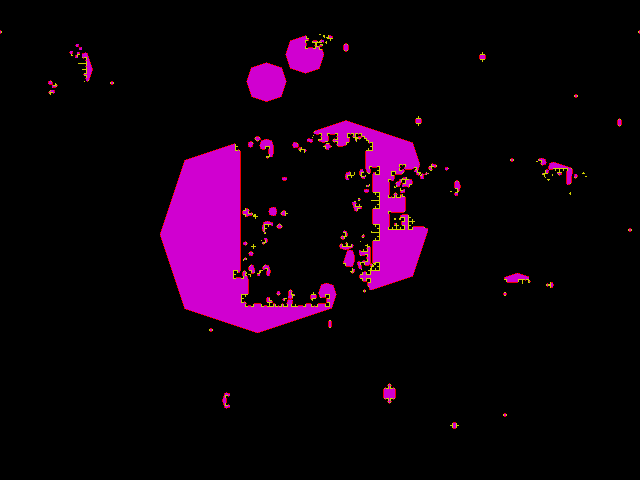}
\includegraphics[width=1in]{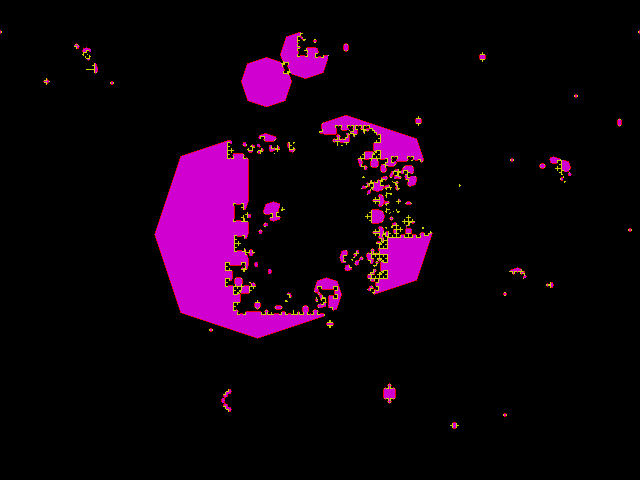}
\end{center}

\begin{center}
\includegraphics[width=1in]{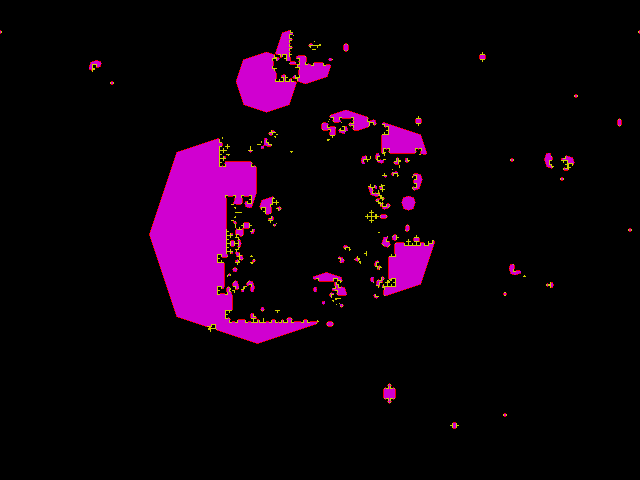}
\includegraphics[width=1in]{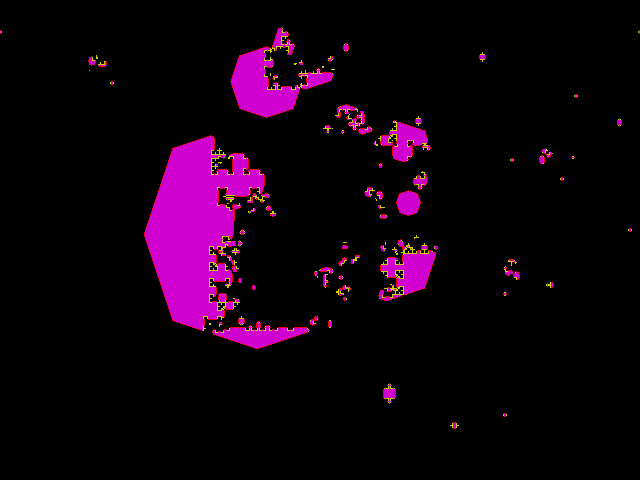}
\includegraphics[width=1in]{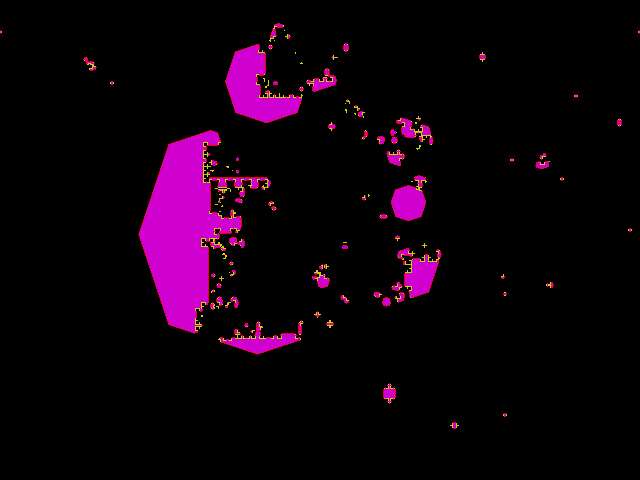}
\includegraphics[width=1in]{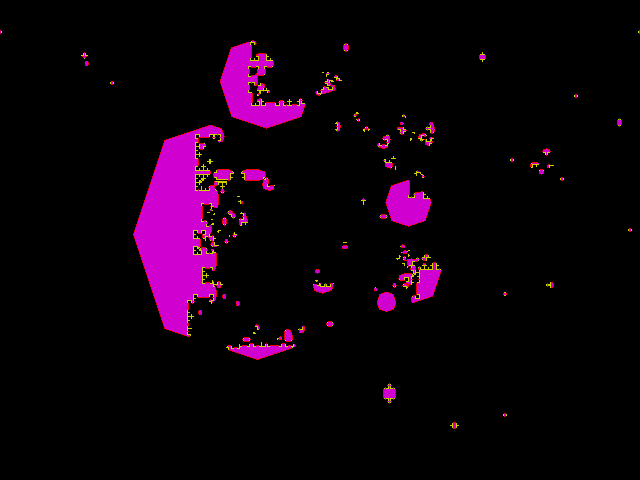}
\end{center}

with a movie available here:

\tt https://dl.dropboxusercontent.com/u/37273687/octo\_02.mp4 \rm

Any nonlinear spatial game can be formulated as a cellular automata, where the future of any particular cell is a function of that cell and its neighborhood. Due to the two-step nature of spatial games, this involves neighbors of neighbors. For a rectangular grid this would be a $5\times 5$ square. The set of all possible ranking matrices for a rectangular grid is $18!$, which would be a minuscule proportion of of the possible cellular automata. Since the nonlinear spatial games rely on counts, the cellular automata are already invariant under symmetry.

Exploring the dynamics of randomly generated rank matrices reveals a remarkable proportion that generate visually interesting dynamics.

An additional example:

\begin{tabular}{c|c|c|c|c|c|c|c|c|c|}
\multicolumn{1}{c}{} & \multicolumn{9}{c}{\hbox{Number of neighbors of type 1}}\\
\multicolumn{1}{c}{} &\multicolumn{1}{c}{0}&\multicolumn{1}{c}{1}&\multicolumn{1}{c}{2}&\multicolumn{1}{c}{3} &\multicolumn{1}{c}{4}&\multicolumn{1}{c}{5}&\multicolumn{1}{c}{6}&\multicolumn{1}{c}{7}&\multicolumn{1}{c}{8}\\
\cline{2-10}
\hbox{ Player type 0:\ } &9&18&4&13&5&1&8&7&14 \\
\cline{2-10}
\hbox{ Player type 1:\ } &3&12&10&17&11&6&16&15&2 \\
\cline{2-10}
\end{tabular}

with movie:

\tt https://dl.dropboxusercontent.com/u/37273687/cellz\_02.mp4 \rm

A hexagonal example:

\begin{tabular}{c|c|c|c|c|c|c|c|}
\multicolumn{1}{c}{} & \multicolumn{7}{c}{\hbox{Number of neighbors of type 1}}\\
\multicolumn{1}{c}{} &\multicolumn{1}{c}{0}&\multicolumn{1}{c}{1}&\multicolumn{1}{c}{2}&\multicolumn{1}{c}{3} &\multicolumn{1}{c}{4}&\multicolumn{1}{c}{5}&\multicolumn{1}{c}{6}\\
\cline{2-8}
\hbox{ Player type 0:\ } &4&13&1&5&10&2&7  \\
\cline{2-8}
\hbox{ Player type 1:\ } &9&14&12&11&6&8&3 \\
\cline{2-8}
\end{tabular}

with movie:

\tt https://dl.dropboxusercontent.com/u/37273687/turq\_07.mp4 \rm

\end{document}